\documentclass[10pt]{amsart}
\usepackage{a4}
\usepackage{amsfonts, amssymb}
\usepackage[all]{xy}

\newtheorem{thm}{Theorem}
\newtheorem{cor}{Corollary}

\newtheorem{defn}{Definition}

\numberwithin{equation}{section} \numberwithin{thm}{section}
\numberwithin{lem}{section} \numberwithin{problem}{section}
\numberwithin{cor}{section}

\begin{document}
\title{Visible lattice points and the chromatic zeta function of a graph}
\address{Instituto de Ciencias Matem\'aticas
 (ICMAT) and Departamento de Matem\'aticas,
Universidad Aut\'onoma  de Madrid,
28049 Madrid, Spain}

\begin{abstract}
We study the probability that a random polygon of $k$ vertices in the lattice $\{1,\dots,n\}^s$ does not contain more lattice points than the  $k$ vertices of the polygon. Then we introduce the chromatic zeta function of a graph to generalize this problem to other configurations induced by a given graph $\mathcal H$.
\end{abstract}
\author{Javier Cilleruelo}

\thanks{This work has been supported by MINECO project MTM2014-56350-P  and ICMAT
Severo Ochoa project SEV-2011-0087}

\maketitle

\section{Introduction} Two distinct points $X,Y$ of the s-dimensional integer lattice are said to be mutually visible if the line segment joining them contains no other lattice point.
 We denote this situation by $X\diamond Y$.
It is well known \cite{Ch} that if  $X,Y$ are lattice points taken at random uniformly in $[1,n]^s$ then $\mathbb P(X\diamond Y)\sim \zeta^{-1}(s)$ as $n\to \infty$, where $\zeta(s)$ is the classical  Riemann zeta function. Since  $X\diamond Y$ and $Y\diamond Z$ are independent events, then $\mathbb P(X\diamond Y\diamond Z)\sim  \zeta^{-2}(s)$. 
What about $\mathbb P(X\diamond Y\diamond Z\diamond X)$? In other words, what is the probability that the three edges of a random triangle $X,Y,Z$ contains no other lattice points than their vertices? 

At first sight, we could expect that $\mathbb P(X\diamond Y\diamond Z\diamond X)\sim  \zeta^{-3}(s)$  since  apparently the three events, $X\diamond Y,\ Y\diamond Z,\ Z\diamond X$, are independent events. We  prove that this intuition is not correct. In fact we obtain a more general result.
\begin{thm}\label{main}Let $s, k\ge 2$ positive integers. If $X^1,\dots,X^k$ are lattice points taken uniformly at random  in $[1,n]^s$, we have
$$\lim_{n\to \infty}\mathbb P(X^1\diamond X^2\diamond \cdots \diamond X^k\diamond X^1)=\zeta^{-k}(s)\prod_p\left (1+\frac{(-1)^k}{(p^s-1)^{k-1}}\right ).$$
\end{thm}
The repetition of vertices $X^i$ is allowed in Theorem \ref{main}. However, the probability of these degenerate cases tends to zero as $n\to \infty$, so we could have formulated Theorem \ref{main} saying that $X^1,\dots, X^k$ are distinct lattice points.

It is interesting to note that the value of the limit in Theorem \ref{main} is smaller than $\zeta^{-k}(s)$ when $k$ is odd and greater than $\zeta^{-k}(s)$ when $k$ is even. We do not understand  the reason of this phenomenom. 

The following version of Theorem \ref{main} can be more illustrative. Take a lattice point $X^1$ at random. Then take a random lattice point $X^2$  visible from $X^1$, then take a random lattice point $X^3$ visible from $X^2$ and so on. What is the probability that $X^k$ is visible from $X^1$? Corollary \ref{cor}, which is a trivial consequence of Theorem \ref{main},  answers this question.
\begin{cor}\label{cor}
Let $s, k\ge 2$ positive integers. If $X^1,\dots,X^k$ are lattice points taken uniformly at random in $[1,n]^s$ we have 
$$\lim_{n\to \infty}\mathbb P(X^k\diamond X^1/X^1\diamond X^2\diamond \cdots \diamond X^k)= \zeta^{-1}(s)\prod_p\left (1+\frac{(-1)^k}{(p^s-1)^{k-1}}\right ).   $$
\end{cor}
Again, we see that $\mathbb P(X^k\diamond X^1/X^1\diamond X^2\diamond \cdots \diamond X^k)$ is smaller than $\mathbb P(X^k\diamond X^1)$ when $k$ is odd and greater when $k$ is even. 
Theorem \ref{main} can be extended to more general configurations.
\begin{defn}
Given a graph $\mathcal H$ of order $k$ we say that that a sequence of lattice points $(X^1,\dots,X^k)$ is $\mathcal H$-visible if $X^i \diamond X^j$ whenever $\{i,j\}\in E(\mathcal H)$. 
\end{defn}Our main Theorem is the following.
\begin{thm}\label{main2}
Let $s, k\ge 2$ positive integers and  $\mathcal H$ a graph  of order $k$. If $X^1,\dots,X^k$ are lattice points taken uniformly at random  in $[1,n]^s$ then we have 
$$\lim_{n\to \infty}\mathbb P((X^1,\cdots, X^k)\text{ is } \mathcal H\text{-visible})= \zeta_{\mathcal H}^{-1}(s), $$
where $\zeta_{\mathcal H}(s)$ is the chromatic zeta function of $\mathcal H$ defined by $$\zeta_{\mathcal H}(s)=\prod_p\left (\frac{P_{\mathcal H}(p^s)}{p^{ks}}\right )^{-1}$$ where $P_{\mathcal H}$ is the chromatic polynomial of $\mathcal H$.
\end{thm}

If we consider the linear graph $\mathcal H=L_k$, with chromatic polynomial $P_{L_k}(x)=x(x-1)^{k-1}$, we have that $\zeta_{\mathcal H}(s)=\zeta^{k-1}(s)$
and we recover the classic result \cite{Ch}:
$$\lim_{n\to \infty}\mathbb P(X^1\diamond  \cdots \diamond X^k:\ X^i\in [1,n]^s)=\zeta^{-(k-1)}(s).$$

 Theorem \ref{main} follows from Theorem \ref{main2} by taking the cycle of $k$ vertices, $\mathcal H=C_k$, and observing that $P_{C_k}(x)=(x-1)^{k}+(-1)^k(x-1): $
$$\zeta_{C_k}^{-1}(s)=\prod_p\left (\frac{(p^s-1)^k+(-1)^k(p^s-1)}{p^{ks}}   \right )=\prod_p\left (1-\frac 1{p^s}\right )^k\prod_p\left (1+\frac{(-1)^k}{(p^s-1)^{k-1}} \right ).$$

David Rearick \cite{R} considered a related problem. Given a set  $S_m=\{X^1,\dots,X^m\}$ of $m$ mutually visible lattice points, he studied the probability that a random lattice point in $[1,n]^s$ is visible from all the lattice points of $S_m$. He proved that 
\begin{equation}\label{R}\lim_{n\to \infty}\mathbb P(X\in [1,n]^2:\ X\diamond X^i,\ i=1,\dots,m)=\prod_p\left (1-\frac m{p^s}\right )\end{equation}
if $m<2^s$ and $0$ if $m\ge 2^s$. In particular \eqref{R} implies that if $m<2^s$ and $X^1,\dots,X^{m+1}$ are taken uniformily at random in $[1,n]^s$, then  \begin{eqnarray*}\label{PP}& &\lim_{n\to \infty}\mathbb P(X^1,\dots,X^{m+1} \text{ is } K_{m+1}\text{-visible}/ X^1,\dots,X^{m} \text{ is } K_{m}\text{-visible})=\prod_p\left (1-\frac m{p^s}\right ).\end{eqnarray*}
This result can be obtained easily from Theorem \ref{main2} considering the chromatic polynomials of the complete graphs,
 $P_{K_{m+1}}(x)=x(x-1)\cdots (x-m),\quad P_{K_{m}}(x)=x(x-1)\cdots (x-m+1)$, and observing that
\begin{eqnarray*}& &\frac{\zeta_{K_{m+1}}^{-1}(s)}{\zeta_{K_m}^{-1}(s)}= \prod_{p}\frac{P_{K_{m+1}}(p^s)}{p^{(m+1)s}}\prod_p\frac{p^{ms}}{P_{K_m}(p^s)}=\prod_p\frac{p^s-m}{p^s}=\prod_p\left (1-\frac m{p^s}\right ).\end{eqnarray*} 

%
%
%

\section{Proof of Theorem \ref{main2}}
Given two lattice points $X^i=(x_1^i,\dots,x_s^i)$ and $X^j=(x_1^j,\dots,x_s^j)$ we write $X^i \equiv X^j\pmod p$ if $x_r^i \equiv x_r^j \pmod p$ for all $r=1,\dots,s$. We write $X^i\not \equiv X^j\pmod p$ otherwise.

Given a prime $p$, we say that $(X^1,\dots,X^k)$ is $\mathcal H_p$-visible if $X^i\not \equiv X^j\pmod p$ whenever $\{i,j\}\in E(\mathcal H)$.
 The first observation is that  \begin{equation}\label{iff}(X^1,\dots,X^k)\text{ is }\mathcal H\text{-visible } \iff (X^1,\dots,X^k)\text{ is }\mathcal H_p\text{-visible for any prime }p.\end{equation}
 
For any positive integer $M$ and $n>M$ we have
\begin{eqnarray}\label{ss}& &|\{X^1,\dots,X^k\in [1,n]^s:\ (X_1,\dots,x^k)  \text{ is }\mathcal H_p\text{-visible for any }p\}|\\&= &|\{ X^1,\dots,X^k\in [1,n]^s:\ (X^1,\dots, X^k) \text{ is }\mathcal H_p\text{-visible for any }p\le M\}|+O(|R|),\nonumber
\end{eqnarray}
where
$$R=\{ X^1,\dots,X^k\in [1,n]^s:\ (X^1,\dots, X^k) \text{ is not }\mathcal H_p\text{-visible for some }p>M\}.  $$
We split $R$ in two sets: $R=R_1\cup R_2$. The set $R_1$ contains those $(X^1,\dots,X^k)$ with $X^i=X^j$ for some $i\ne j$ and $R_2$ contains those with all $X^i$ distinct.

Clearly,
\begin{equation}\label{aa}|R_1|\le \binom k2 n^{s(k-1)}.\end{equation}

On the other hand we observe that if $X^i\ne X^j$ then $X^i\not \equiv X^j\pmod p$ for $p\ge n$, so $(X^1,\dots,X^k)$ is always $\mathcal H_p$-visible when $p\ge n$ for those $(X^1,\dots, X^k)$ counted in $R_2$. Indeed, for a fixed $X^i=(x_1^i,\dots,x_s^i)$ the number of $X^j=(x_1^j,\dots, x_s^j)\in [1,n]^s$ such that $X^j\equiv X^i\pmod p$ is $(n/p+O(1))^s\ll n^s/p^s$ for $p<n$. Thus,
 \begin{eqnarray*}\label{M}
|R_2| & \le &\sum_{M<p< n}|\{\text{distinct }X^1,\dots,X^k\in [1,n]^s:\ (X^1,\dots, X^k) \text{ is not }\mathcal H_p\text{-visible  }\}|\nonumber \\ &\le &\sum_{M<p<n}|\{\text{distinct }X^1,\dots,X^k\in [1,n]^s:\ X^i\equiv X^j\pmod p \text{ for some }i\ne j\}\}|\nonumber\\
&\le &\sum_{M<p<n}\binom k2|\{\text{distinct }X^1,\dots,X^k\in [1,n]^s:\ X^1\equiv X^2\pmod p \}\}|\nonumber\\& \ll& \sum_{M<p<n}\frac{n^{ks}}{p^s} \end{eqnarray*}
and we get the upper bound
\begin{equation}\label{M}  |R_2|\ll n^{ks}M^{1-s}.\end{equation}

By \eqref{ss}, \eqref{aa} and \eqref{M} we have
 \begin{eqnarray}\label{sss}& &|\{ X^1,\dots,X^k\in [1,n]^s:\ (X^1,\dots, X^k) \text{ is }\mathcal H\text{-visible }\}|\\&=&|\{X^1,\dots,X^k\in [1,n]^s:\ (X^1,\dots, X^k) \text{ is }\mathcal H_p\text{-visible for any }p\le M\}|\nonumber \\ & &+O(n^{s(k-1)})+O(n^{ks}M^{1-s}).\nonumber
\end{eqnarray}

The next step is to estimate the quantity
\begin{equation}\label{Hp}|\{X^1,\dots,X^k\in [1,n]^s:\ (X^1,\dots, X^k) \text{ is }\mathcal H_p\text{-visible for any }p\le M\}|.\end{equation}  

A good coloration of a labeled graph $\mathcal H$ is an assigment of colours to the vertices such that two adjacents vertices do not share the same colour. The polynomal chromatic $P_{\mathcal H}(x)$ counts the number of good colorations of $\mathcal H$ using $x$ colours.

For each $p$ we assign to each vertex $X=(x_1,\dots,x_s)$ the $p$-colour $c_p(X)$ defined as the only vector $c_p(X)\in[0,p-1]^s$  such that $c_p(X)\equiv X\pmod p$.

 We observe that $(X^1\dots,X^k)\text{ is } \mathcal H_p\text{-visible  }$ if and only if there exists a good $p$-coloration $C_p=(c_p^1,\dots,c_p^k)$ of $\mathcal H$ such that $ (c_p(X^1),\dots,c_p(X^k))=C_p.$ 
  
Thus, $(X^1\dots,X^k)\text{ is } \mathcal H_p\text{-visible  }$ for any $p\le M$ if and only if there exists a sequence of good colorations $(C_p)_{p\le M}$  such that $ (c_p(X^1),\dots,c_p(X^k))=C_p$ for all $p\le M$.

Since for each prime $p$ there are $p^s$ colours, the number of good $p$-colorations of $\mathcal H$ is $P_{\mathcal H}(p^s)$, where $P_{\mathcal H}$ is the chromatic polynomial of $\mathcal H$. Therefore, the number of sequences of good colorations $(C_p)_{p\le M}$ is 
\begin{equation}
\label{numbercolour}
\prod_{p\le M}P_{\mathcal H}(p^s).\end{equation}
Thus we have 
\begin{eqnarray}\label{pp1}
& &|\{X^1,\dots,X^k\in [1,n]^s:\ (X^1,\dots, X^k) \text{ is }\mathcal H_p\text{-visible for any }p\le M\}|\\ &=&\sum^*|\{X^1,\dots,X^k\in [1,n]^s:\ (c_p(X^1),\dots,c_p(X^k))=C_p,\ p\le M\}|\nonumber
\end{eqnarray}
where the sum $\sum^*$ is extended over all sequences of good colorations $(C_p)_{p\le M}$  of the graph $\mathcal H$.


Given a sequence of colorations $(C_p)_{p\le M}=(c_p^1,\dots,c_p^k)_{p\le M}$ we have that
\begin{eqnarray}\label{pp2}& &|\{X^1,\dots,X^k\in [1,n]^s:\ c_p(X^i)=c_p^i,\ i=1,\dots,k,\text{ for all } p\le M\}|\\ &=&\prod_{i=1}^k|\{X\in [1,n]^s:\ c_p(X)=c_p^i,\text{ for all } p\le M\}|.\nonumber\end{eqnarray}

Given the vectors $c_p^i=(c_{p1}^i,\dots,c_{ps}^i),\ p\le M$, the lattice points $X=(x_1,\dots,x_s)$ with $c_p(X)=c_p^i$ for all $p\le M$ will be those such that the congruences $x_r\equiv c_{pr}^i\pmod p,\ p\le M$ hold for any $r=1,\dots,s$. By the Chinese Remainder Theorem these congruences are equivalent, for each $r=1,\dots,s$, to the congruence $x_r\equiv a_r\pmod{\prod_{p\le M}p}$ for some $a_r$. The number of $x_r\le n$ satisfying each congruence is $\frac n{\prod_{p\le M}p  }+O(1),$
so the number of $X\in [1,n]^s$ with $c_p(X)=c_p^i$ for all $p\le M$ is
$$\left (\frac n{\prod_{p\le M}p  }+O(1)\right )^s.$$
Since this estimate does not depend on the values of $c_p^i$ we have
\begin{equation}\label{pp3}\prod_{i=1}^k|\{X\in [1,n]^s:\ c_p(X)=c_p^i\text{ for all }p\le M\}|=\left (\frac n{\prod_{p\le M}p  }+O(1)\right )^{sk}.\end{equation}

Summing up, as consequence of \eqref{pp1}, \eqref{pp2}, \eqref{pp3} and \eqref{numbercolour} we obtain

\begin{eqnarray*}& &|\{X^1,\dots,X^k\in [1,n]^s:\ (X^1,\dots, X^k) \text{ is }\mathcal H_p\text{-visible for any }p\le M\}|\\&=&\left (\frac n{\prod_{p\le M}p  }+O(1)\right )^{sk}\times |\{\text{sequences of good colorations }(c_p^1,\dots,c_p^k),\ p\le M\}| \\ &=&\left (\frac n{\prod_{p\le M}p  }+O(1)\right )^{sk}\prod_{p\le M}P_{\mathcal H}(p^s)=n^{sk}\left (\prod_{p\le M}\frac{P_{\mathcal H}(p^s)}{p^{sk}}\right )\left (1+O\left(\frac{\prod_{p\le M}p}{n} \right )\right )^{sk}.\end{eqnarray*}
In terms of probability we have proved that
\begin{eqnarray*}&&\mathbb P\left (\{X^1,\dots,X^k\in [1,n]^s:\ (X^1,\dots, X^k) \text{ is }\mathcal H_p\text{-visible for any }p\le M\}\right )\\ &=& \left (\prod_{p\le M}\frac{P_{\mathcal H}(p^s)}{p^{sk}}\right ) \left (1+O\left(\frac{\prod_{p\le M}p}{n} \right )\right )^{sk}.  \end{eqnarray*}
Using \eqref{sss}  we have that
\begin{eqnarray*} & &\mathbb P\left (\{X^1,\dots,X^k\in [1,n]^s:\ (X^1,\dots, X^k) \text{ is }\mathcal H\text{-visible }\}\right )\\&=& \left (\prod_{p\le M}\frac{P_{\mathcal H}(p^s)}{p^{sk}}\right ) \left (1+O\left(\frac{\prod_{p\le M}p}{n} \right )\right )^{sk}+O(n^{-s})+O(M^{1-s}).\end{eqnarray*} Taking the limit as $n\to \infty$ we get
$$ \lim_{n\to \infty}\mathbb P\left (\{X^1,\dots,X^k\in [1,n]^s:\ ()X^1,\dots, X^k) \text{ is }\mathcal H\text{-visible }\}\right )= \prod_{p\le M}\frac{P_{\mathcal H}(p^s)}{p^{sk}}+O(M^{1-s}).$$
Finally, taking the limit as $M\to \infty$ we have
$$ \lim_{n\to \infty}\mathbb P\left (\{X^1,\dots,X^k\in [1,n]^s:\ (X^1,\dots, X^k) \text{ is }\mathcal H\text{-visible }\}\right )= \prod_{p}\frac{P_{\mathcal H}(p^s)}{p^{sk}}=\zeta_{\mathcal H}^{-1}(s).$$

\end{document}